\magnification=\magstephalf
\input amstex
\documentstyle{amsppt}
\pagewidth{6.5truein}
\pageheight{8.9truein}
\ifx\refstyle\undefinedZQA\else\refstyle{C}\fi

\define\R{{\bold R}}
\define\X{{\Cal X}}
\define\Z{{\bold Z}}
\define\nullset{\varnothing}
\define\QED{ \null\nobreak\hfill$\blacksquare$}
\define\QEd{ \null\nobreak\hfill$\square$}
\define\procl#1{\smallskip{\it {#1}. }}

\define\kk{{\bar k}}
\define\mm{{\bar m}}

\define\axis{\ell}
\define\orbit{{\Cal O}}
\define\ident{1}

\define\ri{\therosteritem}
\define\scong{\preceq}
\define\supcong{\succeq}
\define\DF#1/{Dougherty and~Foreman~\cite{\DoughertyForeman#1}}
\define\setgraph{{\Cal G}}

\define\Adams{1}
\define\BanachTarski{2}
\define\DekkerDSS{3}
\define\DekkerFP{4}
\define\Dougherty{5}
\define\DoughertyForeman{6}
\define\Hausdorff{7}
\define\Kuratowski{8}
\define\Magnus{9}
\define\Oxtoby{10}
\define\Robinson{11}
\define\Wagon{12}

\topmatter
\title Solutions to congruences using sets with the property of
Baire\endtitle
\rightheadtext{Solutions to congruences using Baire sets}
\author Randall Dougherty\endauthor
\affil Ohio State University\endaffil
\date February 19, 1999 \enddate
\address Department of Mathematics, Ohio State University,
Columbus, OH 43210\endaddress
\email rld\@math.ohio-state.edu \endemail

\abstract
Hausdorff's paradoxical decomposition of a sphere with countably many
points removed (the main precursor of the Banach-Tarski paradox)
actually produced a partition of this set into three pieces $A,B,C$ such
that $A$~is congruent to~$B$ (i.e., there is an isometry of the set
which sends $A$ to~$B$), $B$~is congruent to~$C$, and $A$~is congruent
to $B \cup C$.  While refining the Banach-Tarski paradox, R.~Robinson
characterized the systems of congruences like this which could be
realized by partitions of the sphere with rotations witnessing the
congruences: the only nontrivial restriction is that the system should
not require any set to be congruent to its complement.  Later, Adams
showed that this restriction can be removed if one allows arbitrary
isometries of the sphere to witness the congruences.

The purpose of this paper is to characterize those systems of
congruences which can be satisfied by partitions of the sphere or
related spaces into sets with the property of Baire.  A paper of
Dougherty and~Foreman gives a proof that the Banach-Tarski paradox can
be achieved using such sets, and gives versions of this result using
open sets and related results about partitions of spaces into congruent
sets.  The same method is used here; it turns out that only one
additional restriction on a system of congruences is needed to make it
solvable using subsets of the sphere with the property of Baire (or
solvable with open sets if one allows meager exceptions to the
congruences and the covering of the space) with free rotations witnessing the
congruences. Actually, the result applies to any complete metric space
acted on in a sufficiently free way by a free group of homeomorphisms.
We also characterize the systems solvable on
the sphere using sets with the property of Baire but allowing all
isometries. \endabstract

\endtopmatter
\document

\head 1. Introduction and definitions \endhead

The basic form of the Banach-Tarski paradox can be stated as follows:
The two-dimensional sphere~$S^2$ can be partitioned into finitely many
pieces $A_1,A_2,\dots,A_n,\allowbreak B_1,B_2,\dots,B_m$ with the
property that the sets~$A_i$ can be rearranged by rigid motions
(rotations) so as to cover the entire sphere, and so can the sets~$B_j$.
This contradicts standard intuitions concerning measure or area; the
sets $A_i$ and~$B_j$ cannot all be measurable with respect to the
standard rotation-invariant probability measure on~$S^2$.

This result of Banach and Tarski~\cite{\BanachTarski}
was based on earlier work of Hausdorff~\cite{\Hausdorff, p.~469}
who proved that the free product of cyclic groups $\Z_2$ and~$\Z_3$ can
be embedded in the rotation group~$SO_3$ of~$S^2$.  Using this,
Hausdorff showed that there is a countable set~$D$ such that
$S^2\setminus D$ can be partitioned into three sets $A,B,C$ such that
$A$ is congruent to~$B$ (i.e., there is a rotation $\rho$ such that
$\rho(A) = B$), $B$~is congruent to~$C$, and $C$ is congruent to $A \cup
B$.  This also is counterintuitive, and the sets $A$,~$B$, and~$C$
cannot be measurable with respect to the standard isometry-invariant
probability measure on~$S^2$. (In fact, there is no rotation-invariant
finitely additive probability measure on~$S^2$ which assigns a measure
to these sets.)

Later, R.~Robinson~\cite{\Robinson} refined the Banach-Tarski construction,
and characterized the systems of congruences for which one can partition~$S^2$
(without a countable exceptional set) into pieces satisfying the
congruences.  In order to state Robinson's results precisely,
we need some definitions.

Fix a positive integer~$r$.  A {\it congruence\/} is specified by two
subsets $L$ and~$R$ of $\{1,2,\dots,r\}$, and is written formally as
$\bigcup_{k\in L} A_k \cong \bigcup_{k\in R} A_k$, where $A_1,A_2,\dots,A_r$
are variables.  The congruence is {\it proper\/} if both $L$ and~$R$
are nonempty proper subsets of $\{1,\dots,r\}$.  Now suppose $G$~is a
group acting on a set~$X$, and a system of congruences is given by pairs
$L_i,R_i\subseteq\{1,\dots,r\}$ for $i\le m$; a {\it solution\/} to the
system of congruences in~$X$ is a sequence of sets $A_k\subseteq X$
($k\le r$) which are pairwise disjoint and have union~$X$, such that, for each
$i\le m$, there is $\sigma_i\in G$ such that $\sigma_i(\bigcup_{k\in
L_i} A_k) = \bigcup_{k\in R_i} A_k$ (i.e., $\sigma_i$~witnesses
congruence number~$i$).  It is clear that only proper congruences are
useful here, and that, if $\sigma_i$~witnesses the congruence given by
$L_i$ and~$R_i$, then $\sigma_i$~also witnesses the complementary
congruence given by $L_i^c$ and~$R_i^c$, where $S^c = \{1,\dots,r\}
\setminus S$.  Also, congruence is transitive: if $\sigma$~witnesses $A
\cong B$ and $\tau$~witnesses $B \cong C$, then $\tau \circ \sigma$
witnesses $A \cong C$.  A system of congruences is called {\it weak\/}
if, among all congruences which can be deduced from the system by taking
complements and applying transitivity, there is no congruence of the
form $\bigcup_{k\in L} A_k \cong \bigcup_{k\in L^c} A_k$ which requires
some set to be congruent to its complement.

It is easy to see that, if a system of congruences has a solution in~$S^2$
with rotations witnessing the congruences, then the system must be weak:
any rotation has fixed points, and hence cannot witness that a set is
congruent to its complement.  Robinson showed that the converse is true:
any weak system of congruences has a solution in~$S^2$ with rotations
witnessing the congruences.  Dekker~\cite{\DekkerDSS} gave
the following abstract
form of Robinson's result:  If $G$~is a free group on more than one
generator which acts locally commutatively on a set~$X$, then any
weak system of congruences has a solution in~$X$ with elements of~$G$
witnessing the congruences.  (An action of a group~$G$ is called
{\it locally commutative\/} if any two elements of~$G$ with a common fixed
point commute; clearly the rotation group on~$S^2$ has this property.
Variants or corollaries of Hausdorff's embedding of a free product
of $\Z_2$~and~$\Z_3$ into~$SO_3$ show that free groups on any finite
or countable number of generators can be embedded into~$SO_3$.)
Also, Adams showed that, if one allows arbitrary isometries of~$S^2$
rather than just rotations to witness the congruences, then any
system of proper congruences has a solution in~$S^2$.

The preceding information is from Wagon~\cite{\Wagon}, which is an
excellent reference on the Banach-Tarski paradox and related work.

In general, a weak system of congruences need not have a solution using
{\sl measurable\/} subsets of~$S^2$.  For example, consider the system
$A_2 \cong A_2 \cup A_3 \cup A_4$, $A_4 \cong A_1 \cup A_2 \cup A_4$,
which Robinson used to get a minimal Banach-Tarski decomposition of~$S^2$.
If these congruences are satisfied by measurable subsets of~$S^2$,
then, since $S^2$~has finite measure, the first congruence forces
$A_3$ and~$A_4$ to have measure~$0$, and the second forces $A_1$ and~$A_2$
to have measure~$0$, so the four sets together cannot cover~$S^2$.
(The same argument applies to Hausdorff's system of congruences, but
this is not weak.)  Another example is the system $A_1 \cong A_2 \cong
A_3 \cong A_4 \cong A_5$, $A_1 \cup A_2 \cong A_1 \cup A_3 \cup A_4$;
any measurable solution to the first part of this would have to give
measure~$1/5$ to each of the sets, making the final congruence
impossible.

If one considers solutions using sets with the property of Baire instead
of measurable sets, then the situation is quite different; the arguments
of the preceding paragraph do not apply.  It was shown in \DF/
that the Banach-Tarski paradox can
be carried out using pieces with the property of Baire. In the present
paper, the methods of \DF/
will be used to characterize those systems of congruences which have
solutions in~$S^2$ under rotations using sets with the property of
Baire.  The result applies more generally, to show that suitable systems
of congruences have solutions using sets with the property of Baire in
any Polish space (complete separable metric space)~$\X$ with a
nonabelian free group of homeomorphisms of~$\X$ which acts locally
commutatively on~$\X$ and freely (without fixed points) on a comeager
subset of~$\X$.  To specify which systems are `suitable' requires
further definitions.

We will call a system of congruences {\it nonredundant\/} if no
congruence in the system can be deduced from the other congruences in
the system by complementation and transitivity as above, and there is no
identity congruence $A \cong A$ in the system.

Next, say that $A$ is {\it subcongruent} to $B$ ($A \scong B$) if $A$ is
congruent to a subset of $B$.  From a given system of congruences, one
can deduce subcongruences by the following rules:  if $A \subseteq B$,
then $A \scong B$; if $A \scong B$ and $B \scong C$, then $A \scong C$;
and, if $A \cong B$ is in the given system, then $A \scong B$, $B \scong
A$, $A^c \scong B^c$, and $B^c \scong A^c$ (where $A^c$~is the
complement of~$A$).  We will call the system of congruences {\it
consistent\/} if there do not exist sets $L,R \subseteq \{1,2,\dots,r\}$
with $R$~a proper subset of~$L$ such that one can deduce $\bigcup_{k \in
L} A_k \scong \bigcup_{k \in R} A_k$ from the system. For example, the
systems used by Hausdorff and Robinson as above are not consistent, but
the other example system above is consistent (the only subcongruences
deducible from it where the left side is a union of more sets than the
right side are $A_1 \cup A_3 \cup A_4 \scong A_1 \cup A_2$ and $A_3 \cup
A_4 \cup A_5 \scong A_2 \cup A_5$).

The main result of this paper is that, if a given system of $m$~congruences
is weak and consistent, and if $\X$~is a Polish space on
which a free group~$G$ of homeomorphisms with $m$~generators acts
locally commutatively everywhere and freely on a comeager set, then the
system of congruences has a solution on~$\X$
using sets with the property of Baire; furthermore, if the system
is nonredundant, then one can use a specified list of $m$~free
generators of~$G$ to serve as the witnesses for the congruences.  (This
latter condition holds for the Robinson-Dekker construction, without
any extra assumption.)  The conditions of weakness and consistency are
necessary for the case of~$S^2$ with a free group of rotations, at least
if we require the sets to be nonmeager; without this requirement, a
system has a solution in~$S^2$ using free rotations if and only if it
has a subsystem (obtained by deleting zero or more of the sets
$A_1,\dots,A_r$ from all congruences) which is weak and consistent.
Also, requiring that a redundant congruence be witnessed by a free
rotation can make a system unsolvable on~$S^2$.

As in \DF/, the results here
concerning sets with the property of Baire are obtained by combining
known results about arbitrary sets with new results about open sets. In
most cases, one cannot expect to get actual solutions to systems of
congruences using open sets; in particular, a connected space cannot be
nontrivially partitioned into open sets at all.  We will therefore allow
meager exceptional sets when trying to satisfy congruences using open
sets.  This leads to the following definitions:  Suppose $G$~is a group
of homeomorphisms of a space~$\X$.  Two sets $A,B\subseteq \X$ will be
called {\it quasi-disjoint\/} if their intersection is meager.  (Of
course, quasi-disjoint open sets in a Polish space are actually
disjoint.)  Sets $A$ and~$B$ are {\it quasi-congruent}, as witnessed by
$\sigma \in G$, if $\sigma(A)$ differs from~$B$ by a meager set.  A {\it
quasi-solution\/} to a system of congruences $\bigcup_{k\in L_i} A_k
\cong \bigcup_{k\in R_i} A_k$ is a sequence of sets $A_k\subseteq \X$
($k\le r$) which are
pairwise quasi-disjoint and whose union is a comeager subset
of~$\X$, such that, for each $i\le m$, there is $\sigma_i\in G$ which
witnesses that $\bigcup_{k\in L_i} A_k$ is quasi-congruent to
$\bigcup_{k\in R_i} A_k$.

The remainder of this paper is as follows.  In section~2, it will be
shown that, if $G$~is a suitable free group of homeomorphisms of a
Polish space $\X$, then any weak consistent system of congruences has a
quasi-solution in~$\X$ using nonempty open sets (with specified free
generators of~$G$ witnessing the congruences, if the system is
nonredundant).  This result is then combined with the results of
Robinson and Dekker to produce solutions (not just quasi-solutions) to
any weak consistent system using sets with the property of Baire.
Section~3 shows the necessity of weakness, consistency, and
nonredundancy.  Section~4 gives the proof that, if one allows arbitrary
isometries of~$S^2$ as witnesses, then any consistent system of
congruences has a quasi-solution using nonempty open sets, and a
solution using nonmeager sets with the property of Baire.

In a later paper~\cite{\Dougherty}, we
consider the problem of finding open sets which actually satisfy
congruences rather than quasi-congruences (but still are only required
to cover a dense subset of the space, rather than all of it).

We will use the symbol~$\circ$ or simple juxtaposition to denote
a group operation, interchangeably.  All group actions will
be written on the left.
For standard basic facts about free groups, such as the unique
expression of any element as a reduced word in the generators and the
fact that any nonidentity element has infinite order, see any text on
combinatorial group theory, such as Magnus, Karrass, and
Solitar~\cite{\Magnus}.  More advanced facts will be referred to
specifically as needed.

\head 2. Positive results \endhead

\proclaim{Theorem 2.1} Suppose $\X$~is a Polish space and $G$~is a
countable group of homeomorphisms of~$\X$ which acts freely on a
comeager subset of~$\X$, and which has a subgroup which is free on
$m$~generators ($m \ge 1$). Suppose that a system of $m$~congruences is
specified by pairs $(L_i,R_i)$ ($1 \le i \le m$) of subsets of\/
$\{1,2,\dots,r\}$; also suppose that this system is weak and consistent.
Then there is a sequence of nonempty open sets $A_k \subseteq\X$ ($k \le
r$) which is a quasi-solution to the system.  Furthermore, if the system
is nonredundant, and elements $f_i$ ($1\le i\le m$) of~$G$ are free
generators for a free subgroup of~$G$, then there is a sequence of
nonempty open sets~$A_k$ as above such that, for each $i \le m$,
$f_i$~witnesses that\/ $\bigcup_{k\in L_i} A_k$ is quasi-congruent to\/
$\bigcup_{k\in R_i} A_k$. \endproclaim

\demo{Proof} First note that, if one congruence in a system is deducible
from the other congruences, then one can delete that one congruence to
get a smaller system, and any quasi-solution to the smaller system will
be a quasi-solution to the original system.  (The same holds for
solutions.)  By iterating this, one can reduce the original system to a
nonredundant system with the same quasi-solutions.  Therefore, it will
suffice to prove only the second part of the theorem.  We may assume
that $G$~is the free group generated by the elements~$f_i$.

We will follow the method of \DF/, but with a few differences.  One
difference is that we will concentrate on the points to be excluded
from the sets~$A_k$, and not construct the sets~$A_k$ themselves until
the excluded sets are complete.  (The reason for this is that it is
easier to work with congruences between intersections than with
congruences between unions; we can actually make intersections
congruent, rather than quasi-congruent.) We will construct sets $B_k$
($1 \le k \le r$) with the following properties: $\bigcap_{k = 1}^r B_k
= \nullset$; the sets $\bigcap_{k' \ne k} B_{k'}$ for $k \le r$ are all
nonempty, and their union is dense in~$\X$; and, for each $i \le m$,
$f_i(\bigcap_{k \in L_i} B_k) = \bigcap_{k \in R_i} B_k$. Once we have
these sets, we can define $A_k$ to be $\bigcap_{k' \ne k} B_{k'}$; then
the sets~$A_k$ will be as desired.  (The intersection of any two
sets~$A_k$ will be $\bigcap_{k' = 1}^r B_{k'}
= \nullset$.  For any $L \subseteq \{1,2,\dots,r\}$, the
set $\bigcap_{k \in L} B_{k}$ includes~$A_{k'}$ for ${k'} \notin L$
and is disjoint from~$A_{k'}$ for ${k'} \in L$, so $\bigcup_{k \in L} A_k$
differs from the complement of $\bigcap_{k \in L} B_{k}$ by a meager set;
hence, a congruence between $\bigcap_{k \in L} B_{k}$ and
$\bigcap_{k \in R} B_{k}$ yields a quasi-congruence between
$\bigcup_{k \in L} A_k$ and $\bigcup_{k \in R} A_k$.)

The open sets~$B_k$ will be built in stages:  we will construct open
sets $B_k^0 \subseteq B_k^1 \subseteq B_k^2 \subseteq \dotso$ for $k
\le r$ and then let $B_k = \bigcup_{n=0}^\infty B_k^n$.  The sets~$B_k^n$
will satisfy the following properties, to be maintained as
induction hypotheses:
\roster
\item[2] $\bigcap_{k=1}^r B_k^n = \nullset$.
\item For each $i \le m$, $f_i(\bigcap_{k \in L_i} B_k^n) = \bigcap_{k
\in R_i} B_k^n$ and $f_i(\bigcap_{k \in L_i^c} B_k^n) = \bigcap_{k \in
R_i^c} B_k^n$.
\item For any $x\in \X$, the set of $y\in \X$ which are connected to~$x$
by a chain of active links is finite.
\endroster
(There is no property~\ri1; this numbering is used for compatibility
with \DF/.) Of course, we must define the terms used in~\ri4:

\procl{Definition} Two points $x$ and~$x'$ are {\it linked}, or there is
a {\it link} from~$x$ to~$x'$, if $x' = f_i(x)$ or $x = f_i(x')$ for
some $i \le m$.  Points $x$ and~$x'$ are {\it connected by a chain of
links} if there are points $x_0,x_1,\dots,x_J$ with $x_0 = x$ and $x_J =
x'$ such that there is a link from $x_{j-1}$ to $x_j$ for each $j \le
J$. A link from~$x$ to~$x'$ is {\it active\/} (for the sets~$B_k^n$) if
there is a point in one or more of the sets~$B_k^n$ which is connected
to~$x$ or to~$x'$ by a chain of at most $2^r$~links.

\smallskip
Note that adding one new point to a set~$B_k^{n+1}$ activates only a
finite number of new links, although the finite number is very large.

Let $B_k^0 = \nullset$ for all~$k$; clearly this makes \ri2--\ri4 true
for $n = 0$.  Let $\langle Z_n \colon n = 0,1,2,\dotsc \rangle$ be a
listing of the nonempty sets in some base for the topology of~$\X$; we
may assume that $\X$~itself occurs at least $r$~times in the list. We
must show how to get from~$B_k^n$ to~$B_k^{n+1}$, preserving properties
\ri2--\ri4, so that, for a given nonempty open set $Z=Z_n$, one of the
sets $\bigcap_{k' \ne \kk}B_{k'}^{n+1}$ for $\kk \le r$ will meet~$Z$.
The $t$\snug'th time that $Z = \X$ ($t \le r$), we will set $\kk = t$ in
order to ensure that $\bigcap_{k' \ne t}B_{k'}$ will be nonempty.  Once
this is accomplished for all~$n$, the resulting sets~$B_k$ will have the
desired properties.

So suppose we are given $B_k^n$ ($k \le r$) and $Z = Z_n$.  The first
step is to find a point $x_0 \in Z$ to be put into all but one of the
sets~$B_k^{n+1}$.  In fact, we will find $x_0$ in~$Z'$, where $Z'$ is~$Z$
unless $Z = \X$ and this is one of the first $r$~occurrences of~$\X$
in the list of open sets; in the latter case, if this is the
$t$\snug'th occurrence of~$\X$, then let $Z'$ be the interior of the
complement of~$B_t^n$.  (To see that this set is nonempty, look at a
$G$\snug-orbit on which $G$~acts freely and which does not meet the boundary
of~$B_t^n$; such orbits form a comeager subset of~$\X$.  By~\ri4, some
point in this orbit is not in any of the sets~$B_k^n$, and hence must
be in the interior of the complement of~$B_t^n$.) So $x_0$~will be in~$Z$
in any case.

Let $D$ be the complement of a ($G$\snug-invariant) comeager set on which
$G$~acts freely, and let $D'$ be the union of the images under the elements
of~$G$ of the boundaries of the sets~$B_k^n$; then $D \cup D'$ is
meager.  Let $x_0$ be any point in $Z' \setminus (D \cup D')$.  By~\ri2,
we can find $\kk \le r$ such that $x_0 \notin B_\kk^n$.  In the case
where $Z = \X$ for the $t$\snug'th time ($t \le r$), we have $x_0 \notin
B_t^n$ by the definition of~$Z'$, so we may set $\kk = t$.  We will
ensure that $x_0 \in B_k^{n+1}$ for all $k \ne \kk$; this will take care
of the current density requirement (or the current nonemptiness
requirement).

As in \DF/, we will construct sets~$\hat B_k$ by adding finitely many
points to the sets~$B_k^n$; $\hat B_k$~will be defined to be $B_k^n \cup
\{g(x_0) \colon g \in\nobreak T_k\}$ for some $T_k \subseteq G$.
However, we will
describe the construction a little differently; instead of giving
inductive clauses to define the sets~$T_k$, we will define a set $M_g
\subseteq \{1,2,\dots,r\}$ for each $g \in G$ and then let $T_k = \{g
\in\nobreak G \colon k \in\nobreak M_g\}$.

We define~$M_g$ recursively, based on the reduced form of $g \in G$ in
terms of the generators~$f_i$.  If $g$~is the identity of~$G$, then let
$M_g = \{k \le\nobreak r \colon k \ne\nobreak
\kk\}$.  Otherwise, we can write~$g$
uniquely as $f_i \circ g'$ or $f_i^{-1} \circ g'$ where $g'$~has a
shorter reduced form than $g$~does, and hence $M_{g'}$ is already
defined.  Let $M^+_{g'} = M_{g'} \cup \{k\colon g'(x_0) \in\nobreak B_k^n\}$.
If $M_{g'} = \nullset$, let $M_g = \nullset$.  If $M_{g'} \ne \nullset$
and $g = f_i \circ g'$, then define~$M_g$ as follows: if $L_i \subseteq
M^+_{g'}$, let $M_g = R_i$; if $L_i^c \subseteq M^+_{g'}$, let $M_g =
R_i^c$; otherwise, let $M_g = \nullset$.  (If both $L_i$ and~$L_i^c$ are
subsets of~$M^+_{g'}$, make some arbitrary definition such as $M_g =
\{1,2,\dots,r\}$; we will see in the next paragraph that this case
cannot occur.)  If $M_{g'} \ne \nullset$ and $g = f_i^{-1} \circ g'$,
then define~$M_g$ in the same way, but with $L_i$ and~$R_i$
interchanged.

First, we show by induction on $g \in G$ that $M_g^+ \ne
\{1,2,\dots,r\}$. If $g$~is the identity, then $\kk \notin M_g^+$ by the
definition of~$x_0$. Otherwise, we have $g = f_i \circ g'$ or $g =
f_i^{-1} \circ g'$ for some simpler~$g'$.  If $M_g = \nullset$, then
$M_g^+ = \{k\colon g(x_0) \in\nobreak B_k^n\} \ne \{1,2,\dots,r\}$ by~\ri2.
If $M_g \ne \nullset$, $g = f_i \circ g'$, and $L_i \subseteq M^+_{g'}$,
then $L_i^c \not\subseteq M^+_{g'}$ by the induction hypothesis, so
$L_i^c \not\subseteq \{k\colon g'(x_0) \in\nobreak B_k^n\}$, so $R_i^c
\not\subseteq \{k\colon g(x_0) \in\nobreak B_k^n\}$ by~\ri3, so $M^+_g = R_i
\cup \{k\colon g(x_0) \in\nobreak B_k^n\} \ne \{1,2,\dots,r\}$.  The remaining
cases are handled the same way.

We can now check that, for any $g$ and~$g'$ in~$G$ and $i \le r$, if $g
= f_i \circ g'$, then $L_i \subseteq M^+_{g'}$ iff $R_i \subseteq
M^+_g$, and $L_i^c \subseteq M^+_{g'}$ iff $R_i^c \subseteq M^+_g$.
First, suppose the reduced form of~$g'$ does not have~$f_i^{-1}$ as its
leftmost term; then $M_g$~is defined from~$M_{g'}$ as above.  If
$M_{g'} = \nullset$ and hence $M_g = \nullset$, then these two
equivalences follow directly from~\ri3; so suppose $M_{g'} \ne
\nullset$.  Now the two left-to-right implications are immediate.  For
the first right-to-left implication, if $L_i \not\subseteq M^+_{g'}$,
then $R_i \cap M_g = \nullset$ by definition of~$M_g$, while $R_i
\not\subseteq \{k\colon g(x_0) \in\nobreak B_k^n\}$ because otherwise
\ri3~would
give $L_i \subseteq \{k\colon g'(x_0) \in\nobreak B_k^n\} \subseteq
M^+_{g'}$, so $R_i \not\subseteq M^+_g$.  The other implication is
proved in the same way.  This completes the case where $g'$~does not
have~$f_i^{-1}$ as its leftmost term.  If $g'$~does have~$f_i^{-1}$ as
its leftmost term, then $g$~does not have~$f_i$ as its leftmost term,
so we can write $g' = f_i^{-1} \circ g$ and proceed as above.

We are now ready to prove \ri2--\ri4 for the sets~$\hat B_k$.  The
definitions of $T_k$ and~$\hat B_k$ (and the fact that $G$~acts freely
on the orbit of~$x_0$) easily imply that $\{k \colon g(x_0) \in\nobreak \hat
B_k\} = M^+_g$ for all $g \in G$, while $\{k \colon x \in\nobreak \hat B_k\} =
\{k \colon x \in\nobreak B_k^n\}$ if $x$~is not in the $G$\snug-orbit of~$x_0$.
Therefore, properties \ri2 and~\ri3 for~$\hat B_k$ follow from the same
properties for~$B_k^n$ and the above facts about~$M^+_g$.

To prove property~\ri4 for the sets~$\hat B_k$, we will need the
following claims, which are the part of this proof where all of the
restrictions on the system of congruences are needed.

Define a labeled directed graph~$\setgraph$ from the system of congruences
as follows.  The vertices of~$\setgraph$ are the nonempty proper subsets
of $\{1,2,\dots,r\}$.  If $S$~is such a subset, and $L_i \subseteq S$,
then $\setgraph$~has an edge from~$S$ to~$R_i$ labeled~$f_i$.
If $L_i^c \subseteq S$, then $\setgraph$~has an edge from~$S$ to~$R_i^c$
labeled~$f_i$.  Similarly, if $R_i \subseteq S$, then
$\setgraph$~has an edge from~$S$ to~$L_i$ labeled~$f_i^{-1}$;
if $R_i^c \subseteq S$, then
$\setgraph$~has an edge from~$S$ to~$L_i^c$ labeled~$f_i^{-1}$.

The digraph~$\setgraph$ has cycles of length~$2$ connecting pairs
$(L_i,R_i)$ or $(L_i^c,R_i^c)$; each such cycle consists of an
$f_i$\snug-edge and an $f_i^{-1}$-edge.  Call the edges in these
$2$\snug-cycles (the edges which come from actual congruences rather
than subcongruences) {\sl good} edges, and call all other edges (e.g., an
$f_i$\snug-edge from a proper superset of~$L_i$ to~$R_i$) {\sl bad} edges.

\procl{Claim 1} No cycle in~$\setgraph$ contains a bad edge.

\procl{Proof} Suppose the edges $e_1,e_2,\dots,e_J$ form a cycle.
Let $N_0,N_1,\dots,N_J$
be the vertices of this cycle (with $N_J = N_0$), so that
$e_j$~is an edge from~$N_{j-1}$ to~$N_j$.
For each~$j$, $e_j$~has a label,
which is either $f_{i_j}$ or~$f_{i_j}^{-1}$ for some~$i_j$.
We will
abuse notation by writing $M \scong N$ for $M,N \subseteq
\{1,2,\dots,r\}$ to mean that the subcongruence $\bigcup_{k \in M} A_k
\scong \bigcup_{k \in N} A_k$ is deducible from the given system of
congruences by the usual rules.

For each~$j$ such that $0 < j \le J$, define a set~$N^-_j$ as follows.
First, suppose $e_j$~is labeled~$f_{i_j}$.
Then either $L_{i_j} \subseteq N_{j-1}$ and
$N_j = R_{i_j}$, or $L_{i_j}^c \subseteq N_{j-1}$ and $N_j =
R_{i_j}^c$; let $N^-_j$ be $L_{i_j}$ in the former case and $L_{i_j}^c$
in the latter.  Similarly, if $e_j$~is labeled~$f_{i_j}^{-1}$,
let $N^-_j$ be $R_{i_j}$
or~$R_{i_j}^c$, depending on whether $N_j$ is $L_{i_j}$ or~$L_{i_j}^c$.
In any case, we have $N_{j-1} \supseteq N^-_j$, and
congruence number~$i_j$ relates either the sets $N^-_j$ and~$N_j$
or their complements.  Therefore, $N_{j-1} \supcong N^-_j \supcong N_j$
for all~$j$.  Since $\supcong$~is transitive and $N_J = N_0$, we have
$N^-_j \supcong N_{j-1}$ for $0 < j \le J$.  Since the system of congruences
is consistent, $N^-_j$ cannot be a proper subset of~$N_{j-1}$, so we must
have $N^-_j = N_{j-1}$ for $0 < j \le J$; this means that all of the edges
are good.\QEd\par\smallskip

Now construct the undirected graph~$\setgraph_0$ by treating each pair
of oppositely-directed good edges in~$\setgraph$ as a single
undirected edge between $L_i$ and~$R_i$ or between $L_i^c$
and~$R_i^c$.

\procl{Claim 2} The undirected graph~$\setgraph_0$ is acyclic (i.e.,
its connected components are trees).

\procl{Proof}  Note that sets $N,N'$ are in the same component
of~$\setgraph_0$ if and only if the congruence $\bigcup \{A_k \colon k
\in\nobreak N\} \cong \bigcup \{A_k \colon k \in\nobreak N'\}$
follows from the given system of congruences.  In particular,
$N$~and~$N^c$ cannot be in the same component of~$\setgraph_0$,
because the system is weak.
Also, note that congruence number~$i$
gives rise to two edges of~$\setgraph_0$,
one between $L_i$ and~$R_i$ and one between $L_i^c$ and~$R_i^c$;
these edges must be in different components of~$\setgraph_0$.

Now suppose we have a nontrivial cycle in~$\setgraph_0$; by taking
a minimal such cycle, we may ensure that there are no repeated
edges in the cycle.  Let one of the edges in the cycle
be an edge from~$L$ to~$R$, coming from congruence number~$i$.
Then the rest of this cycle cannot use this edge and cannot use
the other edge coming from congruence number~$i$ (since this is
not even in the same component), so it consists entirely of edges
coming from other congruences.  But the rest of the cycle gives a path
from~$L$ to~$R$, so $\bigcup \{A_k \colon k
\in\nobreak L\} \cong \bigcup \{A_k \colon k \in\nobreak R\}$
is deducible from the system without using congruence number~$i$.
So congruence number~$i$ is deducible from the others, contradicting
the assumption that the system is nonredundant.\QEd\par\smallskip

Using these two claims, we get:

\procl{Claim 3} Every path of length~$2^r$ in the
digraph~$\setgraph$ contains a pair of consecutive edges
with labels $f_i$ and~$f_i^{-1}$, or vice versa, for some~$i$.

\procl{Proof}  Suppose we have a path of length~$2^r$ in~$\setgraph$.
Since there are fewer than~$2^r$ vertices in~$\setgraph$, some
vertex must be visited more than once, so we get a nontrivial subpath
which starts and ends at the same vertex (i.e., a cycle).  By
Claim~1, this subpath consists entirely of good edges, so
it induces a corresponding path in the graph~$\setgraph_0$
which also starts and ends at the same place.  By Claim~2,
this latter path cannot be a nontrivial cycle, so it must
double back on itself (use the same edge twice in succession);
hence, the original path uses both edges of a pair of oppositely-directed
good edges successively, which gives the desired conclusion.
\QEd\par\smallskip

Now, for any $g \in G$, $x_0$~is connected to~$g(x_0)$ by a chain
of links, and this chain can be read off from the reduced form of~$g$.
In order to prove~\ri4 for the sets~$\hat B_k$, it will suffice to show
that, if $M_g \ne \nullset$, then either all of the links in this chain are
active for the sets~$B_k^n$, or the chain has fewer than~$2^r$ links; once
we know this, \ri4 for~$B_k^n$ implies that there are only finitely
many points~$g(x_0)$ such that $M_g \ne \nullset$ (equivalently, since
$G$~acts freely on the orbit of~$x_0$, the set of~$g$ such that
$M_g \ne \nullset$ is finite), so only finitely many new links are
activated when $B_k^n$~is enlarged to~$\hat B_k$, so \ri4 for~$B_k^n$
implies \ri4 for~$\hat B_k$.

So suppose $M_g \ne \nullset$ and the above chain has at least $2^r$~links.
Then $M_h \ne \nullset$ for all of the intermediate points~$h(x_0)$
on the chain. It must now be true that, given any $2^r$~consecutive
links in the chain, at least one of the $2^r+1$ endpoints of
these links is in one of the sets~$B_k^n$, because otherwise the sets~$M_h$
at these $2^r+1$ endpoints would give a counterexample to
Claim~3. (If none of these points~$h(x_0)$ is in any of the sets~$B_k^n$,
then we
always have $M^+_h = M_h$.  Now, if $h$ and $h' = \rho \circ h$ are
final subwords of the reduced word for~$g$,
where $\rho$~is $f_i$ or~$f_i^{-1}$,
then the way in which $M_{h'}$~is computed from~$M_h$ shows that
there is an edge in~$\setgraph$ from~$M_h$ to~$M_{h'}$ labeled~$\rho$.
The resulting path of length~$2^r$ cannot
include consecutive edges labeled $f_i$ and~$f_i^{-1}$ or vice versa
because we are working with the reduced form of~$g$.)
It follows that all~$2^r$
of the links are active for~$B_k^n$; since this was an arbitrary
subchain of the chain, all of the links in the chain are active for~$B_k^n$.
This completes the proof of~\ri4 for~$\hat B_k$.

Now that we have \ri2--\ri4 for~$\hat B_k$, we can enlarge these sets
to get open sets.  Let $S$ be the set of $g \in G$ such that $x_0$~is
connected to~$g(x_0)$ by a chain of links which are active for the
sets~$\hat B_k$, and let $S'$ be the set of $g'\in G$ such that $g'(x_0)$~is
connected to~$g(x_0)$ for some $g \in S$ by a chain of at most
$2^r+1$ links.  Then $T_k \subseteq S$ for all~$k$, $S \subseteq S'$,
and $S$ and~$S'$ are finite by~\ri4.  Let $U_0$ be an
open neighborhood of~$x_0$
so small that the images~$g(U_0)$ for $g \in S'$ are pairwise disjoint and
each of them is either included in or disjoint from each of the sets~$B_k^n$.
(This is possible because, by the choice of~$x_0$, no point
in~$S'$ is on the boundary of any of the sets~$B_k^n$.)  Now let
$B_k^{n+1} = B_k^n \cup \bigcup \{g(U_0) \colon g \in\nobreak T_k \}$
for each~$k$;
we must see that these sets satisfy properties \ri2--\ri4.

From the definition of~$B_k^{n+1}$ and the disjointness of the sets~$g(U_0)$
for $g \in S'$, the following two statements follow easily: If
$x\in g(U_0)$ for some $g\in S'$, then $x\in B_k^{n+1}$ if and only if
$g(x_0)\in \hat B_k$. If $x\in\X$ is not in any of the sets~$g(U_0)$ for
$g\in S$, then $x\in B_k^{n+1}$ if and only if $x\in B_k^n$.

We can now prove \ri2--\ri4 for~$B_k^{n+1}$.

\ri2: If a point~$x$ is in one of the neighborhoods~$g(U_0)$ where $g\in
S$, then $g(x_0)\notin \bigcap_{k=1}^r \hat B_k$ by \ri2 for~$\hat B_k$,
so $x\notin \bigcap_{k=1}^r B_k^{n+1}$; if $x$~is not in one of these
neighborhoods, then $x\notin \bigcap_{k=1}^r B_k^n$ by \ri2 for~$B_k^n$,
so $x\notin \bigcap_{k=1}^r B_k^{n+1}$.

\ri3: We prove $f_i(\bigcap_{k \in L_i} B_k^{n+1}) \subseteq \bigcap_{k
\in R_i} B_k^{n+1}$; the other parts are similar.  Suppose $x \in
\bigcap_{k \in L_i} B_k^{n+1}$.  If $x \in g(U_0)$ for some $g \in S$,
then $g(x_0) \in \bigcap_{k \in L_i} \hat B_k$, so $f_i(g(x_0)) \in
\bigcap_{k \in R_i} \hat B_k$ by \ri3 for~$\hat B_k$; but $f_i \circ g
\in S'$ and $f_i(x) \in f_i(g(U_0))$, so $f_i(x) \in  \bigcap_{k \in
R_i} B_k^{n+1}$.  If $x$~is not in~$g(U_0)$ for any $g \in S$, then $x
\in \bigcap_{k \in L_i} B_k^n$, so $f_i(x) \in \bigcap_{k \in R_i}
B_k^n$ by \ri3 for~$B_k^n$.

\ri4: Let $w$ be any point of~$\X$, and consider the set of all points
connected to~$w$ by a path of links which are active for the
sets~$B_k^{n+1}$. If this set contains no point which is in~$g(U_0)$ for any
$g\in S$, then all of the links connecting the set were in fact active
for $B_k^n$. (Note: If the link from~$x$ to~$x'$ is
activated by~$x''$, because there is a chain of at most~$2^r$ links
connecting~$x''$ to~$x$ or to~$x'$, then all of the links in this chain
are also activated by~$x''$.) Hence, the set is finite by \ri4
for~$B_k^n$.  So suppose $y\in g(U_0)$ is connected by active links to~$w$,
and $g\in S$.  A point is connected to~$w$ if and only if it is
connected to~$y$, so it will suffice to show that only finitely many
points are connected to~$y$.

Suppose $y$~is actively linked to~$y'$, say $y' = f_i(y)$ (the case
$y' = f_i^{-1}(y)$ is similar).  Let $y''$ be a point in one of
the sets $B_k^{n+1}$ such that $y''$~is connected to either
$y$ or~$y'$ by a chain of at most~$2^r$ links.  Then there is an
element~$h$ of~$G$ such that $h(y) = y''$, and the reduced form
of~$h$ in terms of the generators~$f_I$ has length at most $2^r+1$
(and, if it has length $2^r+1$, then the rightmost component is~$f_i$).
Therefore, $h \circ g \in S'$.  We now have $y'' \in h(g(U_0))$,
so, since $y'' \in B_k^{n+1}$, we must have $h(g(x_0)) \in \hat B_k$.
This means that the link from~$g(x_0)$ to $f_i(g(x_0))$ is active for
the sets~$\hat B_k$, so $f_i(y) \in f_i(g(x_0))$ and $f_i \circ g \in S$.

Now this argument can be repeated starting at~$y'$, and so on; the
result is that, for any chain of active (for the sets~$B_k^{n+1}$) links
starting at~$y$, all of the links in the corresponding chain starting
at~$g(x_0)$ are also active (for the sets~$\hat B_k$).  Furthermore, if
$y$~is connected to two different points $y'$ and~$y''$ by such chains of
links, this will give $y'=h'(y)$ and $y''=h''(y)$ for some distinct
elements $h',h''$ of~$G$, and the corresponding points reached
from~$g(x_0)$ will be $h'(g(x_0))$ and $h''(g(x_0))$; since $G$~acts freely
on the orbit of~$x_0$, these two points will also be different.
Therefore, since $g(x_0)$~is connected to only finitely many points, $y$
(and hence~$w$) must be connected to only finitely many points.  This
completes the proof of~\ri4 for the sets~$B_k^{n+1}$.

This completes the induction. \QED\enddemo

One can use this result to give a new proof of Theorem~4.8 from
\DF/:

\proclaim{Corollary 2.2} Suppose $\X$~is a Polish space and $G$~is a
countable group of homeomorphisms of~$\X$ which acts freely on a
comeager subset of~$\X$.  Then, for any $N\ge 3$, if elements~$f_i$
($1\le i\le N$) of~$G$ are free generators for a free subgroup of~$G$,
then there is an open subset~$A$ of~$\X$ such that the sets~$f_i(A)$ are
disjoint and their union is dense in~$\X$.  In fact, if $f_{ij}\in G$
for\/ $1\le i\le j$ and\/ $3\le j\le N$ are free generators for a free
subgroup of~$G$, then there is an open set~$A$ such that, for each~$j$,
the sets $f_{ij}(A)$ for $i\le j$ are disjoint and have dense union.
\endproclaim

\demo{Proof} We will prove the second part; the proof of the first part
can be obtained from this by omitting most of the congruences (in fact,
the first part is essentially a special case of the second).

Let $R$ be the set of sequences $s = \langle s_j \colon 3 \le j \le N
\rangle$ such that $1 \le s_j \le j$ for each~$j$; we will construct a
system of congruences between sets~$A_s$ for $s \in R$.  (Of course, one
can relabel the sets to make the index set $\{1,2,\dots,r\}$, where $r =
|R| = N!/2$.)  The congruences are: $\bigcup\{A_s\colon s(N) =\nobreak
1\} \cong \bigcup\{A_s\colon s(j) =\nobreak i\}$ for each pair $(i,j)\ne
(1,N)$ such that $3 \le j \le N$ and $1 \le i \le j$. The only proper
congruences that can be deduced from this system are those of the form
$\bigcup\{A_s\colon s(j) =\nobreak i\} \cong \bigcup\{A_s\colon s(j')
=\nobreak i'\}$ and their complementary versions; it follows easily
that the system is weak and consistent.  It is also easy to check that
the system is nonredundant. Therefore, by Theorem~2.1, one can find
a quasi-solution to the system using open sets~$A_s$ for $s \in R$,
where the congruence for $(i,j)$ is witnessed by the element $f_{ij}
\circ f_{1N}^{-1}$ of~$G$.  (Since the elements~$f_{ij}$ are free
generators for their subgroup, the elements $f_{ij} \circ f_{1N}^{-1}$
for $(i,j)\ne(1,N)$ are free generators for {\sl their\/} subgroup.)

Now let $A = f_{1N}^{-1}(\bigcup\{A_s\colon s(N) =\nobreak 1\})$.  Then, for
each $(i,j)$, $f_{ij}(A)$ differs from $\bigcup\{A_s\colon s(j) =\nobreak i\}$
by a meager set; it follows that, for each~$j$, the sets $f_{ij}(A)$ for
$i \le j$ are quasi-disjoint and their union is a comeager (hence dense)
subset of~$\X$.  Since quasi-disjoint open sets must actually be disjoint,
we are done.  \QED\enddemo

The trick used here to get congruent rather than quasi-congruent sets is
quite specific; many weak consistent systems
of congruences do not have quasi-solutions in open sets if one actually
requires congruences instead of quasi-congruences.  This will be
explored further in a later paper~\cite{\Dougherty}.

In order to get results concerning sets with the property of Baire,
we will combine the preceding results about open sets with the
Robinson-Dekker results on arbitrary sets, using the following lemma,
which is a variant of Lemma~2.4 from \DF/:

\proclaim{Lemma 2.3} Suppose $\X$~is a Polish space, and $f_1,\dots,f_m$
are homeomorphisms from~$\X$ to~$\X$.  Also suppose that we have a
system of $m$~congruences such that: there is a solution to the system
in~$\X$ with $f_i$~witnessing congruence number~$i$ for $i=1,2,\dots,m$;
and there is a quasi-solution to the system in~$\X$ using nonmeager sets
with the property of Baire so that $f_i$~witnesses congruence number~$i$.
Then there is a solution to the system in~$\X$ using nonmeager
sets with the property of Baire so that $f_i$~witnesses congruence
number~$i$. \endproclaim

\demo{Proof}  Let $G$ be the countable group of homeomorphisms generated by
$f_1,\dots,f_m$.
Suppose the quasi-solution consists of sets~$A'_k$ with
the property of Baire for $1 \le k \le r$, while the solution is given
by sets~$A''_k$, $1 \le k \le r$. Let $D$ be the union of $\X \setminus
\bigcup_{k=1}^r A'_k$, the intersections $A'_k \cap A'_{k'}$ for $k \ne
k'$, and the differences $f_i(\bigcup_{k\in L_i} A'_k) \triangle
\bigcup_{k\in R_i} A'_k$ for $i \le m$. Then $D$~is meager, and so is
the union~$D'$ of all of the images of~$D$ under the elements of~$G$.
Now let $A_k = (A'_k \setminus D') \cup (A''_k \cap D')$.  These sets
have the property of Baire (since $A'_k$ has the property of Baire and
$D'$~and $A''_k \cap D'$ are meager), and they are easily seen to be
disjoint; using the $G$\snug-invariance of~$D'$, it is easy to verify that
$f_i(\bigcup_{k\in L_i} A_k) = \bigcup_{k\in R_i} A_k$ for each~$i$.
Also, since $A'_k$~is nonmeager, $A_k$~is nonmeager.  Therefore, the
sets~$A_k$ are as desired. \QED\enddemo

\proclaim{Theorem 2.4} Suppose $\X$~is a Polish space and $G$~is a
countable group of homeomorphisms of~$\X$ which acts freely on a
comeager subset of~$\X$ and locally commutatively on all of~$\X$, and
which has a subgroup which is free on $m$~generators ($m \ge 1$).
Suppose that a system of $m$~congruences is specified by pairs
$(L_i,R_i)$ ($1 \le i \le m)$ of subsets of $\{1,2,\dots,r\}$; also
suppose that this system is weak and consistent.  Then there is a
sequence of nonmeager sets $A_k \subseteq\X$ ($k \le r$) with the
property of Baire which is a solution to the system. Furthermore, if the
system is nonredundant, and elements~$f_i$ ($1\le i\le m$) of~$G$ are
free generators for a free subgroup of~$G$, then there is a sequence of
sets~$A_k$ as above such that, for each $i \le m$, $f_i$~witnesses that
$\bigcup_{k\in L_i} A_k$ is congruent to $\bigcup_{k\in R_i} A_k$.
\endproclaim

\demo{Proof} The second part follows immediately from Lemma~2.3,
Theorem~2.1, and the results of Robinson and Dekker, while the first follows
from the second as in the proof of Theorem~2.1. \QED\enddemo

We now recall that a free group on two generators has subgroups
which are free on any finite number of generators~\cite{\Magnus,
Problem~1.4.12}.  This allows us to simplify the statement of the
following corollary, which follows from Theorem~2.4 
just as Corollary~2.2 follows from Theorem~2.1:

\proclaim{Corollary 2.5 \rm\cite{\DoughertyForeman, Theorem~5.4}} Suppose
$\X$~is a Polish space and $G$~is a countable group of homeomorphisms
of~$\X$ which acts freely on a comeager subset of~$\X$ and locally
commutatively on all of~$\X$, and which has a subgroup which is free on
more than one generator.  Then, for any $N \ge 3$, $\X$~can be partitioned
into $N$ $G$\snug-congruent pieces with the property of Baire; in fact,
there is a set $A \subseteq \X$ with the property of Baire such that,
for $3\le j\le N$, $\X$~can be partitioned into $j$~pieces congruent
to~$A$. \QED\endproclaim

\head 3. Negative results \endhead

We will now see why the systems of congruences used in Theorems~2.1
and~2.4 must be weak and consistent in order to get the desired
conclusions for all suitable $\X$ and~$G$.  In fact, it suffices to
look at 2.1 only; if one has a solution (or even a quasi-solution) to a
system of congruences using sets~$A_k$ with the property of Baire, and
if $A'_k$~is an open set differing from~$A_k$ by a meager set, then the
sets~$A'_k$ are a quasi-solution to the same system of congruences.  We
will look at the case of the sphere~$S^2$ acted on by a free group of
rotations.  (For other spaces or groups, more systems of congruences
might be solvable.  For example, if $G$~is a free group on countably infinitely
many generators and we put the discrete metric on~$G$, then we get a Polish
space acted on by~$G$ in which every proper system of congruences has a
solution~\cite{\Wagon, Cor.~4.12}, and this solution will automatically
consist of open sets because the space is discrete.)

We first see why consistency is necessary, at least if we want
quasi-solutions involving nonempty open sets. It is easy to verify
that, if the sets~$A_k$ are a quasi-solution to the system, and the
subcongruence $\bigcup_{k \in R} A_k \scong \bigcup_{k \in R'} A_k$ is
deducible from the system, then $\bigcup_{k \in R} A_k$ actually is
quasi-congruent to a subset of $\bigcup_{k \in R'} A_k$.  We now use
the following fact, a slight variant of Proposition~5.5 from \DF/:

\proclaim{Proposition 3.1} If $B$ and~$C$ are quasi-disjoint subsets
of~$S^2$ with the property of Baire such that $B\cup C$ is quasi-congruent
to a subset of~$B$, then $C$~is meager. \endproclaim

\demo{Proof} Let $\sigma$ be an isometry witnessing the
quasi-congruence, and let $B'$ and~$C'$ be the unique regular-open sets
such that $B \triangle B'$ and $C \triangle C'$ are meager. (For the
definition and properties of regular-open sets, see
Oxtoby~\cite{\Oxtoby,~Ch.~4}. The relevant fact is that $B'$~is the
largest open set such that $B \triangle B'$ is meager, and similarly
for~$C'$.)  Since $\sigma(B \cup C) \setminus B$ is meager, $\sigma(B' \cup
C') \setminus B'$ must be meager; but $B'$~is regular-open, so we have
$\sigma(B' \cup C') \subseteq B'$.  Also, since $B$ and~$C$ are
quasi-disjoint, $B' \cap C'$ is a meager open set and hence empty. If
$\lambda$~is the standard rotation-invariant probability measure on~$S^2$,
then $\lambda(B') \ge \lambda(\sigma(B' \cup C')) = \lambda(B'
\cup C') = \lambda(B')+\lambda(C')$, so $\lambda(C') \le 0$, so $C'$~must
be empty, so $C$~is meager. \QED\enddemo

Therefore, if the open sets $A_k \subseteq S^2$ are a quasi-solution
for a system of congruences, and $\bigcup_{k \in R} A_k \scong
\bigcup_{k \in R'} A_k$ is deducible from the system where $R'$~is a
proper subset of~$R$, then $A_k$~must be meager and hence empty for
each $k \in R \setminus R'$.  Hence, if one insists on a solution using
nonempty open sets, then the system of congruences must be consistent.
If it does not matter that some of the sets are empty, then, given a
system of congruences, one should proceed as follows:  Find all
inconsistencies $\bigcup_{k \in R} A_k \scong \bigcup_{k \in R'} A_k$,
$R' \subset R$ deducible from the system, list all indices~$k$
occurring in $R \setminus R'$, and delete the corresponding sets~$A_k$
from the congruences (i.e., delete these indices from $\{1,2,\dots,r\}$
and from all sets $L_i,R_i$ to get new sets $L'_i,R'_i$ defining a new
system of congruences).  This may produce new inconsistencies; if so,
repeat the process, and continue until no inconsistencies remain.  If
nothing is left (all sets~$A_k$ have been declared empty), then the
original system had no quasi-solutions using open subsets of~$S^2$.
Otherwise, the final system is consistent.  If it is also weak, then
the final system, and hence the original system, has solutions for any
$\X$ and~$G$ as in Theorem~2.4; if the final system is not weak, then
we will see below that the final system has no quasi-solutions using
open subsets of~$S^2$ with a free group~$G$ of rotations, and it
follows that the original system also has no quasi-solutions in this
case.

Suppose $G$~is a free group of rotations of~$S^2$, and we restrict
ourselves to congruences which are witnessed in~$G$; we will now see
that only weak systems can have quasi-solutions using open sets in this
case.  (The corresponding statement about solutions using sets with the
property of Baire, or even using arbitrary sets, is trivial because any
rotation has fixed points and hence cannot map a set to its
complement.)  To see this, we use a lemma about open subsets of~$S^2$
which are quasi-invariant under a rotation of infinite order.  (A set~$A$
is said to be {\it quasi-invariant} under a homeomorphism~$f$ if
the symmetric difference $f(A) \triangle A$ is meager.)

\proclaim{Lemma 3.2} If an open subset~$A$ of~$S^2$ is invariant under a
rotation of infinite order around an axis~$\axis$, then $A$~is invariant
under {\it all} rotations around~$\axis$.  The same applies to
quasi-invariance. \endproclaim

\demo{Proof} Let $\sigma$ be a rotation of infinite order around~$\axis$
under which $A$~is invariant.  Fix $x \in A$, and let $C$ be the circle
generated by rotating~$x$ {\sl continuously\/} around~$\axis$; we must
show that $C \subseteq A$.

Let $2\pi\theta$ be the angle through which $\sigma$ rotates~$S^2$;
since $\sigma$~does not have finite order, $\theta$~is irrational.
It follows that the fractional parts of $j\cdot \theta$
for positive integers~$j$ are dense in the interval $(0,1)$;
equivalently, the rotations through positive integer multiples of
$2\pi\theta$ arbitrarily well approximate any rotation around~$\axis$.
In particular, for any $y \in C$, the rotation around~$\axis$ which
takes~$y$ to~$x$ can be approximated by a rotation~$\sigma^j$
so well that $\sigma^j(y)$ is in the open set~$A$; this means
that $y\in \sigma^{-j}(A) = A$.
%We
%now use the standard fact that the fractional parts of $j\cdot \theta$
%for positive integers~$j$ are dense in the interval $(0,1)$.  (Fix
%$\eps > 0$.  By the pigeonhole principle, one can find integers $i <
%i'$ in the range from~$0$ to $1 + 1/\eps$ such that the fractional
%parts of $i \cdot \theta$ and $i' \cdot \theta$ differ by less than
%$\eps$; but the difference must be nonzero because $\theta$ is
%irrational.  Then the fractional parts of $j \cdot \theta$ for $j =
%i'-i,2(i'-i),\dotsc$ sweep through the interval $(0,1)$, coming within
%$\eps$ of every point in this interval.)
%
%Choose $\eps > 0$ so small that all points of~$S^2$ within distance~$\eps$ %of~$x$ are in~$A$.  Any point~$y$ of~$C$ can be obtained by
%rotating~$x$ through some angle $2\pi\phi$, where $0 \le \phi < 1$.
%Choose $j > 0$ such that the fractional part of~$j\theta$ is within
%$\eps/2\pi$ of~$\phi$; then $\sigma^j(x)$ is within~$\eps$ of~$y$, so $y
%\in \sigma^j(A) = A$.
Since $y$~was arbitrary, we have $C \subseteq A$,
as desired.

Now suppose the set~$A$ is just quasi-invariant under~$\sigma$.  Let
$A'$ be the regular-open set that differs from~$A$ by a meager set; then
$\sigma(A')$~is a regular-open set that differs from~$\sigma(A)$ by a
meager set.  But $\sigma(A)$~differs from~$A$ by a meager set, so
$\sigma(A')$~is a regular-open set differing from~$A$ by a meager set,
so it must be equal to~$A'$.  Therefore, $A'$~is invariant under~$\sigma$,
so it is invariant under any rotation~$\tau$ around~$\axis$;
since $A \triangle A'$ and $\tau(A) \triangle \tau(A')$ are meager,
$A$~must be quasi-invariant under~$\tau$. \QED\enddemo

Now, suppose a given system has a quasi-solution in~$S^2$ using a free
group~$G$ of rotations, but is not weak; fix a set $L\subset
\{1,\dots,r\}$ such that the congruence $\bigcup_{k\in L} A_k \cong
\bigcup_{k\in L^c} A_k$ can be deduced from the system.  Then this
quasi-congruence is witnessed by some $\sigma \in G$, which clearly is
not the identity and therefore must be a rotation of infinite order.
But then $\sigma^2$~is also a rotation of infinite order, and
$\bigcup_{k\in L} A_k$ is quasi-invariant under~$\sigma^2$, so it is
quasi-invariant under all rotations around the axis of~$\sigma^2$.  In
particular, $\bigcup_{k\in L} A_k$ is quasi-invariant under~$\sigma$,
so $\bigcup_{k\in L} A_k$ differs from $\bigcup_{k\in L^c} A_k$ by a
meager set, which is impossible because $\bigcup_{k\in L} A_k$ differs
from the complement of $\bigcup_{k\in L^c} A_k$ by a meager set, and
$S^2$~is not meager.  This contradiction shows that the non-weak system
had no quasi-solution after all.

This shows why weakness and consistency are required in Theorem~2.1.
Next, we consider the requirement of nonredundancy for a system of
congruences to be satisfied with specified witnesses to the
congruences.  Even for a simple redundant system such as $A_1 \cong
A_1$, $A_1 \cong A_1$, $A_1 \cong A_2$, $A_1 \cong A_3$, one cannot
arbitrarily specify the witnesses for the congruences: if the first two
congruences are witnessed by rotations of infinite order around
different axes, then Lemma~3.2 implies that $A_1$~must be
quasi-invariant under any rotation around either of these axes, so
$A_1$~must be either empty or comeager.  (By considering the
regular-open set which differs from~$A_1$ by a meager set, we can
reduce this claim to the corresponding claim about invariant sets: if
$A$~is a nonempty open proper subset of~$S^2$, then $A$~cannot be
invariant under all rotations around either of two different axes.
To see this, note that a connected component of~$A$ must have nonempty
boundary.
If $A$~is invariant around an axis, then this boundary consists of one
or two parts, each of which is a point on the axis or a circle obtained
by revolving a point around the axis, so the axis can be reconstructed
given the boundary.) Either of these makes the rest of the system
impossible to satisfy.  With a little more work, one can show that
requiring even a single redundant congruence to be witnessed by a
rotation of infinite order, such as in the system $A_1 \cong A_1$, $A_1
\cong A_2$, $A_1 \cong A_3$, can make a system unsatisfiable.  (It
would require $A_1$ to be a union of spherical disks and annuli with a
common axis, and $A_2$ and~$A_3$ would also have to be such unions but
around different axes; such sets cannot fit together closely enough to
cover a dense subset of~$S^2$.)

Note that, if $G$~is a countable free group of rotations of~$S^2$,
then, since each element of~$G$ other than the identity has only two
fixed points, $G$~acts freely on $S^2 \setminus D$ for some countable
set~$D$.  But $S^2 \setminus D$ is a $G_\delta$~set in~$S^2$ and is
therefore a Polish space itself \cite{\Kuratowski,~\S33~VI}, and the
preceding paragraphs hold for this new space as well.  Hence, even in a
Polish space on which a free group of homeomorphisms acts freely, one
cannot guarantee that a system of congruences has a quasi-solution
using open sets unless the reduced form of that system (after deleting
inconsistencies as above) is weak and consistent.

Another way to modify the space~$S^2$ is as follows:  Let $G$ be a
free group of rotations of~$S^2$ on $\aleph_0$~generators, and let
$z$ be a point of~$S^2$ such that $G$~acts freely on the $G$\snug-orbit
of~$z$; fix another point~$z'$ of~$S^2$ which is neither~$z$ nor
the point opposite~$z$.
For each $g \in G$, consider the
space~$\X_g$ which is the union of~$S^2$ and its tangent ray at~$g(z)$
in the direction of the (shortest) great-circle arc from~$g(z)$ to $g(z')$,
with
the standard Euclidean metric from~$\R^3$.  Note that $h(\X_g) = \X_{h\circ g}$
if we view the rotations as acting on all of~$\R^3$.
Now take a copy of~$\X_g$
for each~$g$, and identify corresponding points of~$S^2$ to get a
space~$\X$.  (Tangent rays that happen to intersect will {\sl not\/} have their common
points identified.  If $x,y \in \X$ are in different tangent rays,
then the distance from~$x$ to~$y$ in~$\X$ is the length of the shortest
path in~$\R^3$ from~$x$ to~$y$ {\sl via a point of~$S^2$}.)  Then $\X$~is
a Polish space, and its group of isometries is precisely~$G$; $G$~acts
locally commutatively on~$\X$ and freely on $\X \setminus D$ for a
countable meager set~$D$,
and the negative results given above for~$S^2$ also apply to~$\X$, so we
can get such results even when using the entire isometry group of a
suitable Polish space.

\head 4. Congruences on the sphere using all isometries \endhead

We now know what congruences have solutions using subsets of~$S^2$ with
the property of Baire and using free rotations; it is natural to ask
what can be done if arbitrary isometries of~$S^2$ are allowed.  The
results in the preceding section concerning consistency apply just as
well for arbitrary isometries, so even here a system must be consistent
(or at least reducible to a consistent system by deletion of some sets)
in order to have such a solution.  However, it turns out that the
restriction of weakness can be removed if we allow arbitrary isometries
to witness the congruences.

As usual, one of the ingredients needed for the proof is a
corresponding result for arbitrary subsets of~$S^2$; this result is due
to Adams~\cite{\Adams} (see also Wagon~\cite{\Wagon,~Theorem~4.16}).
Unfortunately, Adams'
proof cannot be used here; the particular isometries he uses to witness
the congruences cannot be used to get a corresponding result concerning
open sets.  (Adams' construction causes a complementary congruence $A
\cong A^c$ to be witnessed by an isometry~$\tau$ such that $\tau^2$~is
a rotation of infinite order; then $A$~is invariant under~$\tau^2$, and
we saw in the preceding section why this cannot work for open sets.)
We therefore give a revised proof of this result.

\proclaim{Theorem 4.1 \rm(Adams)} Any system of proper congruences has
a solution in~$S^2$, if arbitrary isometries can be used as witnesses
for the congruences. \endproclaim

\demo{Proof}
We first transform the system into an equivalent system having a useful
form.  Call two systems of congruences (on the same index set
$\{1,2,\dots,r\}$) {\it equivalent\/} if any congruence in one can be deduced
from the other, and vice versa; clearly equivalent systems have the same
solutions.  By moving to an equivalent system if necessary, we may
assume that the system is presented with as few congruences as possible
(i.e., there is no equivalent system with fewer congruences).

Now, suppose the system (call it~$S_0$) is not weak.  Let $M_0$ be a
subset of $\{1,2,\dots,r\}$ such that the congruence $\bigcup_{k \in
M_0} A_k \cong \bigcup_{k \in M_0^c} A_k$ is deducible from~$S_0$;
choose~$M_0$ so that this deduction requires as few steps as possible.
Such a deduction gives a sequence $M_0,M_1,\dots,M_n$ of subsets of
$\{1,2,\dots,r\}$ such that $M_n = M_0^c$ and each pair $(M_i,M_{i+1})$
appears as one of the congruences in~$S_0$, perhaps in the complemented
form $(M_i^c,M_{i+1}^c)$.  Because the deduction is minimal, no set
appears more than once in the list $M_0,M_1,\dots,M_n$, and the only
case where both a set and its complement appear in the list is $M_n =
M_0^c$.  But this easily implies that no congruence in~$S_0$ is used
more than once during the deduction: if it were used twice in the same
form, this would require a duplication in the list, while if it were
used once in the given form and once in the complemented form (assuming
these are different), there would be more than one instance of a set
and its complement appearing in the list.  Let $S_1$ be $S_0$~with the
last congruence used in the above deduction (the one between $M_{n-1}$
and~$M_n$, or maybe their complements) deleted, and let $S'_0$ be
$S_1$~together with the congruence $\bigcup_{k \in M_0} A_k \cong \bigcup_{k
\in M_0^c} A_k$.  Since the congruence between $M_0$ and~$M_{n-1}$ is
deducible from~$S_1$, the congruence between $M_n$ and~$M_{n-1}$ is
deducible from~$S'_0$, so $S'_0$~is equivalent to~$S_0$.

Now look at~$S_1$, ignoring the new self-complement congruence.  If
$S_1$~is not weak, one can repeat the above process to change~$S_1$ into
an equivalent system~$S'_1$ with the same number of congruences, where
$S'_1$~is $S_2$~together with another self-complement congruence.
Repeat this process as many times as possible, until we reach a system~$S_j$
which is weak.  Let the congruences in~$S_j$ be given by pairs
$(L_i,R_i)$ for $1 \le r \le \mm$, and let the self-complement
congruences be given as $(L_i,R_i)$ (with $R_i = L_i^c$) for
$\mm+1 \le i \le m$.  (If the original system was weak, then $\mm
= m$.) We have now found a system equivalent to the original system,
with a minimal number~$m$ of congruences (the same number as in~$S_0$),
such that the first~$\mm$ congruences form a weak system and the
remaining $m - \mm$ congruences are between a set and its complement.
Since $m$~is minimal, this system is nonredundant.

The main result of Dekker~\cite{\DekkerFP} states that any
reasonable-sized (continuum or smaller) free product of cyclic groups
can be embedded in the rotation group of~$S^2$.  Therefore, we can
choose rotations $\sigma_i$ ($1 \le i \le \mm$) and $\tau'_i$ ($\mm < i
\le m$) such that each~$\sigma_i$ has infinite order, each~$\tau'_i$ has
order~$4$, and the group~$G'$ generated by all of these rotations is the
free product of the cyclic groups generated by the rotations
individually.  Let $\zeta$ be the antipodal isometry which maps each
point of~$S^2$ to the point opposite it. Let $\tau_i = \zeta \circ
\tau'_i$ for each~$i$, and let $G$ be the group generated by the
isometries $\sigma_i$ and~$\tau_i$.  We will construct a solution to the
revised system of congruences so that $\sigma_i$ ($i \le \mm$) or
$\tau_i$ ($i > \mm$) witnesses congruence number~$i$.

Clearly $\zeta^2$~is the identity on~$S^2$. Since isometries of~$S^2$
preserve oppositeness of points, $\zeta$~commutes with all isometries
of~$S^2$.  Using this, we see that $\tau_i$~has order~$4$ in~$G$, and the
homomorphism from~$G'$ to~$G$ which sends $\sigma_i$ to~$\sigma_i$ and
$\tau'_i$ to~$\tau_i$ (which exists and is unique by the definition of
free products) is in fact an isomorphism.  Therefore, $G$~is also a
free product of $\mm$~copies of~$\Z$ and $m-\mm$ copies of~$\Z_4$.
Also, $G \cap G'$ is a group of index~$2$ in~$G$ and in~$G'$,
consisting of those words in~$G'$ such that the total number of
occurrences of the generators~$\tau'_i$ is even.

We will now show that much of Robinson's work on free groups, as
presented in Chapter~4 of Wagon~\cite{\Wagon}, can be carried out as
well for free products of $\Z$\snug's and $\Z_4$\snug's.  The rest of this
proof will follow the relevant parts of that chapter rather closely.

First, we look at the structure of the group~$G$ (of course, the same
results will apply to~$G'$, which is isomorphic to~$G$).  Any element
of~$G$ has a unique expression as a reduced word.  Here a `word' is a
product (possibly of length~$0$) of elements $\sigma_i^{\pm1}$
and~$\tau_i^{\pm1}$; a word is reduced if there is no occurrence of
$\sigma_i\sigma_i^{-1}$, $\sigma_i^{-1}\sigma_i$, $\tau_i^4$,
or~$\tau_i^{-1}$ (which is equal to~$\tau_i^3$).

Given such a reduced word~$g$, we can express it in the form
$h_1h_2h_3$ with $h_3 = h_1^{-1}$ where $h_1$ and~$h_3$ include as much
of the word~$g$ as possible.
(For a reduced word~$h_1$, the inverse reduced word~$h_1^{-1}$ is
obtained by reversing the word and then replacing $\sigma_i$
with~$\sigma_i^{-1}$, $\sigma_i^{-1}$ with~$\sigma_i$, and maximal
consecutive blocks~$\tau_i^j$ with~$\tau_i^{4-j}$.)
To do this, start by setting $h_1$ and~$h_3$
to be the identity and $h_2$~to be~$g$.  If the current~$h_2$
starts with~$\sigma_i$ and ends with~$\sigma_i^{-1}$, transfer
the~$\sigma_i$ to~$h_1$ and the~$\sigma_i^{-1}$ to~$h_3$; similarly if
$h_2$~starts with~$\sigma_i^{-1}$ and ends with~$\sigma_i$.  If the
current~$h_2$ both starts and ends with one or more copies of~$\tau_i$,
with a total of at least~$4$ such copies (but $h_2$~is not just a power
of~$\tau_i$), then transfer $j$~copies from the start of~$h_2$ to~$h_1$
and $4-j$ copies from the end of~$h_2$ to~$h_3$, where $j$~is $1$,~$2$,
or~$3$, as appropriate.  (If there are more than~$4$ such copies at the
ends of~$h$, so that more than one choice of~$j$ is possible, then use
$j=2$, for a reason to be seen below.)  Repeat until $h_2$~cannot be
reduced further.

When $g$~is expressed as above, it is easy to see that the reduced form
of~$g^n$ is $h_1h_2^nh_3$ for any positive~$n$, unless $h_2$~is of the
form~$\tau_i^k$, in which case the reduced form of~$g^n$ is the null
word (the identity) if $kn$~is divisible by~$4$, $h_1\tau_i^{\text{$kn$
mod $4$}}h_3$ otherwise.  (If $h_2$~is not of the form~$\tau_i^k$, then
the fact that $h_2$~cannot be reduced further as above shows that
$h_2^n$~is a reduced word.)  In fact, the same statement also holds for
negative~$n$, because we used $j=2$ whenever possible in the preceding
paragraph.  Furthermore, if $g^n$~is broken into three pieces as above,
then the three pieces are precisely $h_1$,~$h_2^n$, and~$h_3$.

From these facts, it follows immediately that the only elements of~$G$
of finite order are conjugates of powers~$\tau_i^k$.  In particular, the
only elements of~$G$ of order~$2$ have the form $g\tau_i^2g^{-1}$ for
some $g$ and~$i$.  (A similar statement holds for~$G'$, of course.)

One more fact we will need is that the only abelian subgroups of~$G$
are the cyclic subgroups (so, if two elements of~$G$ commute, then they
are powers of a single element of~$G$).  This follows from the Kurosh
Subgroup Theorem~\cite{\Magnus,~Cor.~4.9.1}.

We now start to work out the analogues for~$G$ of Robinson's results for
free groups.

\proclaim{Lemma 4.2} The group $G$ above can be partitioned into subsets
$A_1,A_2,\dots,A_r$ satisfying the given system of congruences, with
$\sigma_i$ ($i \le \mm$) or $\tau_i$ ($i > \mm$) witnessing congruence
number~$i$ for each $i \le m$.  In fact, for any word~$w$ from~$G$ in
which the total number of occurrences of the generators~$\tau_i$ is even,
there is such a partition of~$G$ which puts~$w$ in the same
set~$A_k$ as the identity element~$\ident$ of~$G$. \endproclaim

\demo{Proof} First, we show that the subsets of $\{1,2,\dots,r\}$ can
be colored with two colors so that: for any set~$L$, $L$ and~$L^c$ have
opposite colors; for any $i \le \mm$, $L_i$ and~$R_i$ have the same
color. To do this, define an equivalence relation on subsets of
$\{1,2,\dots,r\}$ as follows: $L$~is equivalent to~$L'$ iff the
congruence $L \cong L'$ can be deduced from the first~$\mm$ congruences
of the given system.  Clearly, if $L$~is equivalent to~$L'$, then $L^c$~is
equivalent to~$L^{\prime c}$.  Also, $L$~is never equivalent to~$L^c$,
since the first~$\mm$ of the given congruences form a weak
system.  Therefore, the equivalence classes under this relation come in
complementary pairs; if we assign colors to the equivalence classes so
that complementary classes get opposite colors, then the induced
coloring of the subsets of $\{1,2,\dots,r\}$ will be as desired.

We can view the given $m$~congruences as $2m$~formal inclusions:
the equation $\sigma_i(\bigcup_{k \in L_i} A_k) = \bigcup_{k \in R_i} A_k$
can be expressed as the two inclusions $\sigma_i(\bigcup_{k \in L_i}
A_k) \subseteq \bigcup_{k \in R_i} A_k$ and $\sigma_i^{-1}(\bigcup_{k \in
R_i} A_k) \subseteq \bigcup_{k \in L_i} A_k$, and similarly for~$\tau_i$.
We will therefore use the terms `domain of~$\sigma_i$' and `range
of~$\sigma_i$' for the sets $\bigcup_{k \in L_i} A_k$ and $\bigcup_{k
\in R_i} A_k$, respectively, and define `domain of $\tau_i$,' `range of
$\tau_i^{-1}$,' and so on similarly.

Suppose $w = \rho_n\rho_{n-1}\dots\rho_1$, where each~$\rho_k$ is
$\sigma_i^{\pm 1}$ or~$\tau_i$ for some~$i$ (and this is the reduced
form of~$w$).  We will first assign the end segments $\ident, \rho_1,
\rho_2\rho_1, \dots, w$ to suitable sets~$A_k$, and then handle the
remaining elements of~$G$.

First, suppose that, for some $j \le n$, the range of~$\rho_j$ is
neither the domain of~$\rho_{j+1}$ nor the complement of the domain
of~$\rho_{j+1}$ (here we let $\rho_{n+1} = \rho_1$).  Then either the
range of~$\rho_j$ meets both the domain of~$\rho_{j+1}$ and its
complement, or the complement of the range of~$\rho_j$ meets both the
domain of~$\rho_{j+1}$ and its complement; let $S$ be the domain
of~$\rho_j$ in the former case, the complement of this domain in the
latter case.  Assign $\rho_{j-1}\dots\rho_1$ to one of the sets
in~$S$.  Next, assign $\rho_{j-2}\dots\rho_1$ to an appropriate set~$A_k$;
this will be a set in the domain of~$\rho_{j-1}$ if
$\rho_{j-1}\dots\rho_1$ is in the range of~$\rho_{j-1}$, and a set not
in this domain if $\rho_{j-1}\dots\rho_1$ is not in this range.  Repeat
this process to assign all of the shorter end segments of~$w$, down
to~$\ident$, to sets~$A_k$.  Put~$w$ in the same set as~$\ident$, and then
assign $\rho_{n-1}\dots\rho_1$ and so on; continue until only
$\rho_j\dots\rho_1$ remains unassigned.  This word must be assigned to
the range of~$\rho_j$ if $S$~is the domain of~$\rho_j$, the complement
of this range otherwise; it also must be placed in the domain
of~$\rho_{j+1}$ if $\rho_{j+1}\dots\rho_1$ is in the range
of~$\rho_{j+1}$, the complement of this domain otherwise.  By the
definition of~$S$, these two requirements can both be met.  (This must
be reworded slightly in the case $j = n$, but the basic argument
remains the same.)

Now, suppose that the preceding case does not hold; for every~$j$, the
range of~$\rho_j$ is either the domain of~$\rho_{j+1}$ or its
complement.  Let $S_0$ be the domain of~$\rho_1$.  Given a set~$S_{j-1}$
which is either the domain of~$\rho_j$ or the complement of
the domain of~$\rho_j$, let $S_j$ be the range of~$\rho_j$ in the
former case, the complement of the range of~$\rho_j$ in the latter.
Then $S_n$~must be either~$S_0$ or the complement of~$S_0$.  Note that
each~$S_j$ is a union of sets~$A_k$, say $\bigcup_{k \in N_j} A_k$, and
therefore $S_j$ (actually,~$N_j$) has had a color
assigned earlier in the proof of the Lemma.  Furthermore, if $\rho_j$
is~$\sigma_i^{\pm 1}$, then $S_{j-1}$ and~$S_j$ have the same color; if
$\rho_j$ is~$\tau_i$, then $S_{j-1}$ and~$S_j$ have opposite colors.
Since the number of $j$\snug's for which $\rho_j$~is a generator~$\tau_i$ is
even (by hypothesis on~$w$), $S_n$~must have the same color as~$S_0$,
so $S_n$~must be~$S_0$ rather than the complement of~$S_0$.  We now
easily assign each end segment $\rho_j\dots\rho_1$ to one of the sets~$A_k$
included in~$S_j$, making sure to put $\ident$ and~$w$ in the
same set included in~$S_0$; these assignments are compatible with the
required inclusions.

Now that we have assigned the end segments of~$w$ to sets~$A_k$, the
remaining elements~$g$ of~$G$ can be assigned by an easy recursion on
the reduced form of~$g$.  Suppose this reduced form starts with~$\rho$,
where $\rho$~is $\sigma_i^{\pm 1}$ or~$\tau_i$, and $g = \rho g'$ where~$g'$
has already been assigned to a set~$A_k$.  Then, if $g'$~is in the
domain of~$\rho$, assign~$g$ to the range of~$\rho$; if $g'$~is in the
complement of the domain of~$\rho$, assign $g$~to the complement of the
range of~$\rho$.

We must verify that, if $g = \rho g'$ where $\rho$~is $\sigma_i$ or~$\tau_i$,
then $g$~is in the range of~$\rho$ if and only if $g'$~is in
the domain of~$\rho$.  Let $v$ and~$v'$ be the reduced words for $g$
and~$g'$.  If $v = \rho v'$ and $v$~is an end segment of~$w$, then the
way in which the end segments of~$w$ were assigned to sets~$A_k$ gives
the desired result here; the same applies if $v' = \rho^{-1} v$ and
$v'$~is an end segment of~$w$.  If $v = \rho v'$ and $v$~is not an end
segment of~$w$, then we get the desired result from the recursive
definition of the preceding paragraph; this also holds if $v' =
\rho^{-1} v$ and $v'$~is not an end segment of~$w$.  The only remaining
case is when neither $\rho v'$ nor $\rho^{-1} v$ is reduced.  This can
happen only when $\rho$ is~$\tau_i$ and $v' = \tau_i^3 v$.  But then,
by the preceding cases, we have $v \in L_i$ iff $\tau_i v \in L_i^c$
iff $\tau_i^2 v \in L_i$ iff $\tau_i^3 v \in L_i^c$; hence, we have the
desired result in this case as well.  Therefore, the sets~$A_k$ satisfy
the system of congruences. \QED\enddemo

Another useful fact is that, if $w$~is a word in the generators of~$G$
which has an odd number of occurrences of the generators~$\tau_i$ (including
as inverses), then the corresponding isometry of~$S^2$ has no fixed points. 
Let $w'$ be the element of~$G'$ corresponding
to~$w$ (i.e., replace all generators~$\tau_i$ with~$\tau_i'$).  Since
$\tau_i = \zeta \tau_i'$, $\zeta$~commutes with all elements of~$G'$,
$\zeta^2$~is the identity, and the number of occurrences of the
generators~$\tau_i$ in~$w$ is odd, we can compute that
$w = \zeta w'$.
Now $w'$~is a rotation which cannot be of order~$2$, since the only
elements of~$G'$ of order~$2$ are the conjugates of~$\tau_I^{\prime
2}$, which all have even numbers of generators~$\tau'_i$.  But it is
easy to see that the only rotations of~$S^2$ which map some point to
its antipodal point are rotations of order~$2$.  Therefore, $w'$~does
not map any point to its antipodal point, so $w = \zeta w'$ has no
fixed points.

We now resume the proof of Theorem~4.1.  In order to get the desired
partition of~$S^2$, it will suffice to get such a partition for each
$G$\snug-orbit in~$S^2$ and put them together (using the axiom of choice to
choose one such partition for each orbit).  So consider one such
orbit~$\orbit$.  If $G$~acts freely on~$\orbit$, then fixing any element~$x$
of~$\orbit$ determines a bijection $g \mapsto g(x)$ from~$G$ to~$\orbit$
which preserves the action of~$G$, so any partition of~$G$ as
in Lemma~4.2 can be transferred to~$\orbit$, giving a partition of~$\orbit$
with the desired properties.

So suppose $G$~does not act freely on~$\orbit$.  Let $w$ be a
non-identity reduced word of~$G$, as short as possible, such that $w$~has
a fixed point in~$\orbit$; let $x$ be such a fixed point.  Then
$w$~cannot start with~$\sigma_i^{-1}$ and end with~$\sigma_i$, because, if
it did, then the reduced form of $\sigma_i \circ w \circ \sigma_i^{-1}$
would be shorter than~$w$ and would have a fixed point $\sigma_i(x) \in
\orbit$. Similarly, $w$~cannot start with~$\sigma_i$ and end
with~$\sigma_i^{-1}$; and $w$~cannot start with~$\tau_i^j$ and end
with~$\tau_i^{j'}$ where $j + j' \ge 4$, except in the case that $w$~is just
a power of~$\tau_i$.  In fact, if $w$~is not a power of~$\tau_i$, then
we may assume that $w$~does not both start and end with~$\tau_i$; if it
ends with~$\tau_i^j$, then we can replace~$w$ with the reduced form of
$\tau_i^j \circ w \circ \tau_i^{-j}$, which still starts with~$\tau_i$
but ends with something else, and has the fixed point $\tau_i^j(x)$.
We also know that $w$~has an even number of occurrences of
generators~$\tau_i$ (and hence represents a rotation of~$S^2$);
in particular, if
$w$~is a power of some~$\tau_i$, then $w$~must be~$\tau_i^2$.

Let $\rho$ be the leftmost term in the reduced word~$w$ (either
$\sigma_i$, $\sigma_i^{-1}$, or~$\tau_i$ for some~$i$).  Define~$\rho'$
to be~$\sigma_i^{-1}$ if $\rho = \sigma_i$, $\sigma_i$~if $\rho =
\sigma_i^{-1}$, and $\tau_i$~if $\rho = \tau_i$; we have ensured that
$w$~does not end with~$\rho'$, unless $w = \tau_i^2$ for some~$i$.
Therefore, the word~$w^n$ is already in reduced form for positive~$n$,
and the reduced form of~$w^n$ for negative~$n$ does not begin with~$\rho$,
unless $w = \tau_i^2$.

The next thing to show is that the only elements of~$G$ which fix~$x$
are the powers of~$w$.  Suppose $v$~is a nonidentity member of~$G$ such
that $v(x) = x$.  Then $v$~must also have an even number of occurrences
of generators~$\tau_i$,
and is therefore a rotation.  Since the rotation group acts locally
commutatively on~$S^2$, $v$ and~$w$ must commute, so together they
generate an abelian subgroup of~$G$.  We noted earlier that any abelian
subgroup of~$G$ is cyclic, so there must be an element~$u$ of~$G$ such
that $v$ and~$w$ are both powers of~$u$.  We may assume that $w$~is a
{\sl positive} power of~$u$ (replace~$u$ with~$u^{-1}$ if necessary).
If $u$~has finite order, then $u = g^{-1}\tau_i^jg$ for some
$g$,~$i$, and~$j$; since $v$ and~$w$ are nonidentity rotations and are powers
of~$u$, we must have $v = w = g^{-1} \tau_i^2 g$, so $v$~is a power of~$w$.
Now suppose $u$~is of infinite order.  From the general arguments
about the structure of~$G$ given earlier (specifically, the expression
of~$u$ in the form $h_1h_2h_3$ so that the reduced form of~$u^n$ is
$h_1h_2^nh_3$ for any $n > 0$), we see that, if $n > k > 0$, then the
reduced form of~$u^n$ is longer than the reduced form of~$u^k$.  Suppose
$w = u^n$ and $v = u^j$, and let~$k$ be the greatest common divisor of
$n$ and~$j$.  Then $u^k$~can be expressed as a power of~$v$ times a
power of~$w$ (by applying the extended Euclidean algorithm to $n$ and~$j$),
so $u^k(x) = x$.  We clearly have $k \le n$, but we cannot have $k
< n$, since otherwise~$u^k$ would be shorter than~$w$, contradicting the
choice of~$w$ as the shortest possible word with a fixed point in~$\orbit$.
Therefore, $k = n$, so $j$~is divisible by~$n$, so $v$~is a
power of~$w$, as desired.

Using the above, we now show that every element of the orbit~$\orbit$
has a unique expression of the form $g(x)$, where $g$~is an
element of~$G$ whose reduced form does not end in~$w$ and does not end
in~$\rho'$.  (Exception: if $w = \tau_i^2$, then the reduced form of~$g$
is allowed to end in $\rho' = \tau_i$, but not in~$\tau_i^2$.)  Let $y$
be any point in this orbit, and let $v$ be a shortest reduced word in~$G$
such that $v(x) = y$.  Clearly $v$~cannot end in~$w$ (otherwise
the reduced form of
$vw^{-1}$ is shorter).  If $w$~is not of the form~$\tau_i^2$, and
$v$~ends in~$\rho'$, then $vw$~does not end in~$w$, and it does not end
in~$\rho'$ either, since $w$~does not end in~$\rho'$.  (If the entire~$w$
cancels out when $vw$~is transformed to reduced form, then $vw$~has a
shorter reduced form than~$v$, contradicting the choice of~$v$.)
Therefore, we can take $g$ to be either $v$ or~$vw$.  To see that $g$~is
unique, suppose $u$ and~$v$ are distinct and $u(x) = v(x) = y$.  Then
$(u^{-1}v)(x) = x$, so $u^{-1}v$ is a power of~$w$; by interchanging $u$
and~$v$ if necessary, we may ensure that $u^{-1}v$ is a {\sl positive\/}
power of~$w$, say~$w^j$.  Then either $v = uw^j$ ends in~$w$, or there
is some cancellation when $u$~is multiplied by~$w^j$; in the latter
case, $u$~must end in~$\rho'$ (in~$\tau_i^2$ if $w = \tau_i^2$, since in
this case $w^j$~must be~$w$).  This completes the proof that $g$~is
unique.

Now, apply Lemma~4.2 to partition~$G$ into sets $A_1,A_2,\dots,A_r$
satisfying the congruences, so that $w$~is in the same set~$A_k$
as~$\ident$.  This lets us partition~$\orbit$ into sets $B_1,B_2,\dots,B_r$
as follows: for any point $y \in \orbit$, find the unique expression
$g(x)$ for~$y$ as above, and put $y \in B_k$ iff $g \in A_k$.  We must
see that the sets~$B_k$ satisfy the system of congruences.

First, suppose $y \in \orbit$, $i \le \mm$, and $z = \sigma_i(y)$; we
must see that $y \in \bigcup_{k \in L_i} B_k$ if and only if $z \in
\bigcup_{k \in R_i} B_k$. Express $y$ and~$z$ as $g(x)$ and $h(x)$,
where $g$ and~$h$ do not end in~$w$ and (if $w$~is not of the
form~$\tau_I^2$) do not end in~$\rho'$. Then $y \in \bigcup_{k \in L_i} B_k$
iff $g \in \bigcup_{k \in L_i} A_k$, and $z \in \bigcup_{k \in R_i} B_k$
iff $h \in \bigcup_{k \in R_i} A_k$. Also, since the sets~$A_k$ satisfy
the congruences, we have $g \in \bigcup_{k \in L_i} A_k$ iff $\sigma_ig
\in \bigcup_{k \in R_i} A_k$, and $\sigma_i^{-1}h \in \bigcup_{k \in
L_i} A_k$ iff $h \in \bigcup_{k \in R_i} A_k$. Therefore, we are done if
$h = \sigma_ig$. So suppose $h \ne \sigma_ig$.  Then the reduced form of
$\sigma_ig$ must end in~$w$ or in~$\rho'$, while the reduced form of~$g$
does not.  There are only two cases in which this can happen: either
$\rho = \sigma_i$ and $\sigma_ig = w$, or $\rho' = \sigma_i$ and $g =
\ident$.  In the first of these cases we have $h = \ident$, and since
$\ident$ and~$w$ lie in the same set~$A_k$, we have $h \in \bigcup_{k
\in R_i} A_k$ iff $\sigma_ig \in \bigcup_{k \in R_i} A_k$, and this
gives the desired result.  In the second case, we have $\sigma_i^{-1}h =
w$, so $g \in \bigcup_{k \in L_i} A_k$ iff $\sigma_i^{-1}h \in
\bigcup_{k \in L_i} A_k$, and again we are done.

Now suppose $y \in \orbit$, $i > \mm$, and $z = \tau_i(y)$.  Define $g$
and~$h$ as above.  Repeating this argument, we again see that we are
done unless $h \ne \tau_ig$.  Again this happens in only two cases:
either $\rho = \tau_i$ and $\tau_ig = w$, or $\rho' = \tau_i$ and $g =
\ident$ (and $w \ne \tau_i^2$).  These two cases are handled just as
before.

This completes the construction of the desired partition for an
arbitrary orbit of~$S^2$ under~$G$, so we are done.  \QED\enddemo

The corresponding result for open sets is:

\proclaim{Theorem 4.3} Any consistent system of proper congruences
has a quasi-solution in~$S^2$ using nonempty
open sets (and arbitrary isometries).
\endproclaim

\demo{Proof} Revise the system as in the first three paragraphs of the
proof of Theorem~4.1, and define isometries $\sigma_i$,~$\tau'_i$,
$\zeta$, and~$\tau_i$ and groups $G$ and~$G'$ as in the fourth paragraph
of that proof; we will use the same isometries as witnesses here.

The proof will follow that of (the second part of) Theorem~2.1 quite
closely, so we will just give the differences here.  Let $f_i$ be
$\sigma_i$~if $i \le \mm$, $\tau_i$~if $i > \mm$.

The definition of `active link' is changed slightly: a link from~$x$
to~$x'$ will be considered active for the sets~$B_k^n$ if there is a point
in one or more of these sets which is connected to~$x$ or to~$x'$ by a
chain of at most~$2^{r+1}$ (rather than~$2^r$) links.

The next change is at the proof that, if $g = f_i \circ g'$, then $L_i
\subseteq M^+_{g'}$ iff $R_i \subseteq M^+_g$, and $L_i^c \subseteq
M^+_{g'}$ iff $R_i^c \subseteq M^+_g$.  If $f_i$ is~$\sigma_i$, then the
argument is unchanged, but if $f_i$ is~$\tau_i$, then the cases are
slightly different.  If the reduced form of~$g'$ does not start
with~$\tau_i^3$, then $M_g$~is defined from~$M_{g'}$, and we get the desired
result as before.  If the reduced form of~$g'$ does start with~$\tau_i^3$,
so that $g'$~is $\tau_i^3 g$, then $M_{g'}$~is defined from~$M_{\tau_i^2g}$,
which is defined from~$M_{\tau_ig}$, which is defined
from~$M_g$.  Now, using the preceding case and the fact that $R_i =
L_i^c$, we get $L_i \subseteq M^+_{g'}$ iff $L_i^c \subseteq
M^+_{\tau_i^2g}$ iff $L_i \subseteq M^+_{\tau_ig}$ iff $L_i^c \subseteq
M^+_g$; similarly, $L_i^c \subseteq M^+_{g'}$ iff $L_i \subseteq M^+_g$,
as desired.

Next, we must give revised forms of the Claims.
Define a labeled directed graph~$\setgraph$ from the system of congruences
as before, except that the edges labeled $f_i^{-1}$ for $i>\mm$
are omitted.
Again the digraph~$\setgraph$ has cycles of length~$2$ connecting pairs
$(L_i,R_i)$ or $(L_i^c,R_i^c)$, for $i \le \mm$; each such cycle consists of
an $\sigma_i$\snug-edge and an $\sigma_i^{-1}$-edge.  For $i>\mm$
we also get $2$\snug-cycles between $L_i$ and $R_i = L_i^c$;
in this case both edges in the cycle will be labeled~$\tau_i$.
Call the edges in all of these
$2$\snug-cycles {\sl good} edges, and call all other edges {\sl bad} edges.

Since the system is still consistent, the same argument as before gives:

\procl{Claim 1} No cycle in~$\setgraph$ contains a bad edge.
\QEd\par\smallskip

Again construct the undirected graph~$\setgraph_0$ by treating each pair
of oppositely-directed good edges in~$\setgraph$ as a single
undirected edge.

\procl{Claim 2} The undirected graph~$\setgraph_0$ is acyclic; furthermore,
each component of~$\setgraph_0$ contains at most one edge coming
from the self-complement congruences.

\procl{Proof}  Suppose we have a nontrivial cycle in~$\setgraph_0$;
as before, we may assume that this cycle does not use an edge
more than once.  If this cycle includes an edge coming from
congruence number~$i$ for $i > \mm$, then, since congruence number~$i$
produces only one edge of~$\setgraph_0$, the rest of the cycle
must come from the other congruences; as before, this implies
that congruence number~$i$ is deducible from the remaining
congruences, contradicting nonredundancy.  So all of the
edges in the cycle come from the first $\mm$~congruences;
since these congruences form a weak system, we get a contradiction
as in the old Claim~2.

Now, suppose there are two self-complement edges in the same component.
Find a shortest possible path connecting endpoints of two such edges;
this path (possibly of length~0) consists of distinct edges
from the first~$\mm$ congruences.  Say this path connects~$L$
to~$R$, where $L \in \{L_i,L_i^c\}$, $R \in \{L_{i'},L_{i'}^c\}$,
and $i$ and~$i'$ are distinct numbers greater than~$\mm$.
Then there is a nontrivial cycle in~$\setgraph_0$ from~$L$
to~$R$ (the given path) to~$R^c$ (the $\tau_{i'}$\snug-edge)
to~$L^c$ (the complemented and reversed form of the given path)
to~$L$ (the $\tau_i$\snug-edge), contradicting the preceding paragraph.
\QEd\par\smallskip

Using these two claims, we can now get:

\procl{Claim 3} Every path of length~$2^{r+1}$ in the
digraph~$\setgraph$ contains either four consecutive edges
with the same label~$\tau_i$ for some $i>\mm$ or a pair of consecutive edges
with labels $\sigma_i$ and~$\sigma_i^{-1}$, or vice versa, for some~$i\le\mm$.

\procl{Proof}  Suppose we have a path of length~$2^{r+1}$ in~$\setgraph$.
Since there are fewer than~$2^r$ vertices in~$\setgraph$, some
vertex, say~$L$, must be visited at least three times.
Let $p$~be the subpath from~$L$ to~$L$ to~$L$.  By
Claim~1, this subpath consists entirely of good edges, so
it induces a corresponding path~$p_0$ in the graph~$\setgraph_0$
which also goes from~$L$ to~$L$ to~$L$.  By Claim~2,
$p_0$~cannot include a nontrivial cycle, so each
of its two $L$\snug-to-\snug$L$ parts must
double back on itself.
If either doubles back on itself at a $\sigma_i$\snug-edge, then
$p$~has a pair of consecutive edges
with labels $\sigma_i$ and~$\sigma_i^{-1}$, or vice versa, so we are done.
If neither part of~$p_0$ doubles back on itself at a $\sigma_i$\snug-edge,
then they both must double back at a $\tau_i$\snug-edge.  By Claim~2,
there is only one such edge~$e$ in the component of~$\setgraph_0$
containing~$p_0$, and there is a unique path~$q$ in~$\setgraph_0$ from~$L$
to the nearest endpoint of~$e$
which does not double back.  Hence, $p_0$~must consist of
$q$, an even number of traversals of~$e$, $q'$ (the reversal of~$q$),
$q$~again, another even number of traversals of~$e$, and $q'$~again.
If $q$~is non-null, then $p_0$ doubles back on itself at a $\sigma_j$-edge
(the last edge of~$q'$ followed by the first
edge of~$q$), so we are done as before; if $q$~is null, then $p_0$~consists of
at least four consecutive occurrences of~$e$, so $p$~contains
four consecutive $\tau_i$\snug-edges, as desired.
\QEd\par\smallskip

The next step is to show that, for any $g \in G$, if $M_g \ne \nullset$,
then either all of the links in the canonical chain from~$x_0$ to
$g(x_0)$ (i.e., the chain read off from the reduced form of~$g$; note
that there might be other chains from~$x_0$ to~$g(x_0)$, since $G$~is no
longer free) are active for the sets~$B_k^n$, or this chain has fewer
than~$2^{r+1}$ links.  The proof of this works as before (with
$2^r$ replaced by~$2^{r+1}$), so properties \ri2--\ri4 hold for the
sets~$\hat B_k$.

The construction of the sets~$B_k^{n+1}$ goes through as before, with
two minor changes: one must replace $2^r$ with~$2^{r+1}$ throughout, and
one must not assume that there is a unique chain of links connecting two
points in the $G$\snug-orbit of~$x_0$.  The wording of the definition of the
set~$S$ does not need to be changed, but one must note that it refers to
arbitrary chains from~$x_0$ to $g(x_0)$, rather than just the canonical
chain.  Also, in the second paragraph of the proof of~\ri4 for the
sets~$B_k^{n+1}$, one does not necessarily use the {\sl reduced\/} form of
the element of~$h$; instead, one just uses the fact that there is {\sl
some\/} expression of~$h$ as a product of elements~$f_I$ and their
inverses (one is allowed to use~$f_I^{-1}$ even if $I > \mm$) such that
the product has length at most $2^{r+1}+1$, and if its length is equal
to $2^{r+1}+1$, then the rightmost component is~$f_i$.  Everything else
goes through as before.

This completes the induction.\QED\enddemo

Since the same isometries were used in the preceding two proofs to
witness the congruences, Lemma~2.3 now gives:

\proclaim{Theorem 4.4} Any consistent system of proper congruences has a
solution in~$S^2$ using nonmeager
sets with the property of Baire (and arbitrary
isometries). \QED\endproclaim

Actually, the proof of Theorem~4.3 goes through without change if the
involution~$\zeta$ is deleted, so that the group~$G'$ is used instead
of~$G$.  This gives:

\proclaim{Proposition 4.5} Any consistent system of proper congruences
has a quasi-solution in~$S^2$ using nonempty
open sets, with rotations witnessing
the congruences. \QED\endproclaim

However, this result does not lead to a result about sets with the
property of Baire, because there is no corresponding result giving
solutions using arbitrary sets (unless the system is weak).


\Refs

\ref \no \Adams \by J. Adams \paper On decompositions of the sphere
\jour J. London Math. Soc. \vol 29 \yr 1954 \pages 96--99 \endref

\ref \no \BanachTarski \by S. Banach and A. Tarski \paper
Sur la d\'ecomposition des ensembles de points
en parties respectivements congruents \inbook Oeuvres
\bookinfo S. Banach, vol. 1 \publ \'Editions Scientifiques de
Pologne \publaddr Warsaw \yr 1967 \pages 118--148 \endref

\ref \no \DekkerDSS \by Th. Dekker \paper Decompositions of sets and
spaces I \jour Indag. Math. \vol 18 \yr 1956 \pages 581--589 \endref

\ref \no \DekkerFP \by Th. Dekker \paper On free products of cyclic rotation
groups \jour Canad. J. Math. \vol 11 \yr 1959 \pages 67--69 \endref

\ref \no \Dougherty \by R. Dougherty \paper
Open sets satisfying systems of congruences
\toappear \paperinfo arXiv:math.MG/0001010 \endref

\ref \no \DoughertyForeman \by R. Dougherty and M. Foreman \paper
Banach-Tarski decompositions using sets with the property of Baire
\jour J. Amer. Math. Soc. \vol 7 \yr 1994 \pages 75--124 \endref

\ref \no \Hausdorff \by F. Hausdorff \book Grundz\"uge der
Mengenlehre \publ Chelsea \publaddr New York \yr 1949 \endref

\ref \no \Kuratowski \by K. Kuratowski \book Topology \bookinfo vol. 1
\publ Academic Press \publaddr New York \yr 1966 \endref

\ref \no \Magnus \by W. Magnus, A. Karrass, and D. Solitar \book
Combinatorial Group Theory \publ Dover \publaddr New York
\yr 1976 \bookinfo second edition \endref 

\ref \no \Oxtoby \by J. Oxtoby \book Measure and Category \publ
Springer-Verlag \publaddr New York \yr 1980
\bookinfo Graduate Texts in Mathematics 2 \endref

\ref \no \Robinson \by R. Robinson \paper On the decomposition of spheres
\jour Fund. Math. \vol 34 \yr 1947 \pages 246--260 \endref

\ref \no \Wagon \by S. Wagon \book The Banach-Tarski Paradox \publ Cambridge
University Press \publaddr Cambridge \yr 1993 
\bookinfo second edition \endref

\endRefs
\enddocument